%% file: ode.tex
\documentclass[12pt,a4paper]{article}
\usepackage{amsmath}
\usepackage{amsfonts}
\usepackage{bm}
\usepackage{tikz,fullpage}
\usetikzlibrary{backgrounds}
\usetikzlibrary{calc}
\date{}
\begin{document}

\title{Power series with Taylor coefficients of sum-product type and algebraic differential equations /SP-series and the interplay between their resurgence and differential properties}
\author{ Shweta Sharma }
\maketitle

\tableofcontents

\input{ode_1.tex}


\input{ode_2.tex}





\input{ode_7.tex}

\input{ode_8.tex}

\input{ODEbis.tex}

\input{ode_10.tex}


%
\end{document}

%% file: ode_1.tex
\bigskip

\section{Introduction}

\bigskip
The present article is an elaboration with numerical details and some explicit calculations of the chapter 6 of the long paper [ES]. The article is essentially based on the differential properties of the inner generators that occur while handling the SP series (sum product series, a reminder of the definition is given below). The inner generators introduced in [ES] not always but often satisfy linear homogeneous ordinary differential equation. 

The article is structured as follows: First there are just some reminders about basic definitions used in SP series. This serves as a preliminary and a way to settle notation and not the least in the exposition of the problem at hand. Section $2$ is about the two types of differential equation - variable and covariant. We give a schematic description of the existence of such ODEs and recommend [ES] for a more formalistic view. The  next section is about resurgence and connection matrices that appear in our treatment of the problem. Some of the matrices are explicitely written down
for different values of $p$, the parameter in our driving function $f$. In passing we give a list of the $T$ and $\gamma$ polynomials which are related to ODEs of infinite order. The last section talks about the main numerical result of the article i.e. the non existence of differential equations for rational $F$ inputs in a global perspective.

In this series, I'm writing another article [SS] which would be soon finished. It would be essentially involved in treating examples of arbitrary functions $F$ of various types which always fall in the category of the series of the type sum-product (SP) to be treated theoretically as well as numerically, the study of the various singularities that occur, how they could be retrieved using the machinery that is detailed later (using the symmetry between the series and the Taylor coefficient expansion asymptotics) their placement and the resurgence equation and the resurgence that helps us to retrieve the Riemann surface in its complete description using the alien derivatives over the series and the Taylor expansion. We go through some of the analytic properties that we would exploit.

For a polynomial $f$ or a rational function $F$ (resp. a trigonometric polynomial) and Taylor coefficients $ J(n) $ defined only by pure products $\prod $ the series in the $\zeta$-plane would be of hypergeometric (q-hypergeometric) type. Hence, the theory of SP series links two important theories. This problem started with a problem of the study of the series associated to knots[CG3]

\subsection{Definitions and basics}
The sum-product type series are essentially Taylor series of the form $$ j(\zeta) := \sum_{n \geq 0} J(n)\zeta^n $$ whose coefficients are of the sum product (SP) type:
$$ J(n) := \sum_{\epsilon \leq m \leq n} \prod_{\epsilon \leq k \leq m} F(\frac{k}{n}) = \sum_{\epsilon \leq m \leq n} \exp \left(-\sum_{\epsilon \leq k \leq m} f(\frac{k}{n})\right) $$ 
for $ \epsilon \in \{0,1\} $
\\
The two driving functions $F$ and $f$ are related as follows: $ F \equiv \exp(-f) $. The summation always starts at $\epsilon = 0 $ unless $ F(0) \in \{ 0, \infty \}$ (here it starts at $\epsilon = 1$).
$F$ is meromorphic and special care is to be taken in case the function is holomorphic 

%

A few things that need attention are the outer/inner dichotomy and the ingress factor [ES].
The SP series lead to two different types of singularities (generators) as under the alien derivation this generates the resurgence algebra of the SP series. 

The generators that occur in this problem are of the following types:
\begin{itemize}
\item the inner generators, so called as they occur indefinetely under alien derivation but never produce the outer generators. 
\item The outer generators, as they occur one, never twice (they never occur under repeated alien derivation) 
\item The other two generators are less important for the problem, they are the original generators which is the SP series itself and the exceptional generators which do not appear naturally, but which are a very useful tool for the exploration of the Riemann surface, since (i)they produce all inner generators under alien derivation (ii) their base point can be moved around at will and taken as close as one wishes to the base point of any inner generator.
\end{itemize}

\subsection{Definition and calculation of the inner generators.}

\textbf{The \textit{nir}-transform}

The \textit{inner generators} are given by the double integral:  \textit{nir:} $ f \rightarrow h $
\begin{align*}\label{*}
 h(\nu) & = \frac{1}{2 \pi i} \int^{c + i \infty}_{c - i \infty} \exp(n \nu) \frac{dn}{n} \int^{\infty}_{0} \exp^{\#}(-f^{\Uparrow \beta^{*}}(n,\tau))d\tau \\
& = \frac{1}{2 \pi i} \int^{c + i \infty}_{c - i \infty} \exp(n \nu) \frac{dn}{n}\int^{\infty}_{0}\exp(-f_{\kappa}\frac{n^{-\kappa \tau^{ \kappa +1}}}{\kappa +1})\exp_{\#}(-f^{\uparrow \beta^{*}}(n,\tau))d\tau \end{align*}
where $ \exp_{\#}(X)$ (resp. $ \exp^{\#}(X)$) denotes the exponential expanded as a power series of $X$ (resp. minus its leading term). $\beta^{*}_{k} = \frac{B_{k+1}(\tau+\frac{1}{2})}{k+1}$ for the standard case, where $B_k$ is the $k^{\textrm{th}}$ Bernoulli polynomial and $f^{\Uparrow \beta^{*}}(n,\tau)$ and $f^{\uparrow \beta^{*}}(n,\tau)$ are defined as follows:

\begin{align*}
f(x)&= \sum_{k \geq \kappa} f_k x^k \quad (\kappa \geq 1 \, ,\, f_k \not= 0) \\
f^{\Uparrow \beta^{*}} &:= \beta({\partial}_{\tau})f(\frac{\tau}{n}) \quad \textrm{with}\quad {\partial}_{\tau}:=\frac{\partial}{\partial \tau} \\
&:= \sum_{k\geq \kappa} f_k n^{-k} \beta^{*}_{k}(\tau) := f_{\kappa}\frac{n^{-\kappa} {\tau}^{k+1}}{k+1} + f^{\uparrow \beta^{*}}(n,\tau)
\end{align*}
In the integral expression for $h(\nu)$, we first perform the term by term ramified Laplace integration in  $\tau$ followed by a term by term Borel integration in $n$.

The \textit{inner} \textit{generators} also admit a more complicated analytical expression, as a succession of nine successive operations, all of which are more or less elementary (though often highly non-linear) except for the crucial step, the \textit{mir}-transform, which is an integro-differential operator of infinite order (with rational coefficients): see section 4.2 [ES].


\subsection{Definition and calculation of the outer generators.}

\textbf{The \textit{nur}-transform}

The outer generators can be viewed as infinite sums of inner generators and the inner generators (see section 1.2) but they are also capable of a direct construction.

Altogether, we have the choice between three equivalent definitions.

\begin{enumerate}
 \item The first expression: \textit{nur:} $ f \rightarrow h $
\begin{align*}
h(\nu) &=\frac{1}{2 \pi i} \int^{c + i \infty}_{c - i \infty} \exp(n \nu) \frac{dn}{n} \sum^{\infty}_{m=\frac{1}{2} + \mathbb{N}} \exp^{\#}(-f^{\Uparrow \beta^{*} }(n,m)) \\
&= \frac{1}{2 \pi i} \int^{c + i \infty}_{c - i \infty} \exp(n \nu) \frac{dn}{n} \sum^{\infty}_{m=\frac{1}{2} + \mathbb{N}}\exp(-f_{\kappa}\frac{n^{-\kappa} m^{ \kappa +1}}{\kappa+1})\exp_{\#}(-f^{\uparrow \beta^{*}}(n,m))
\end{align*}
is analogous to the double integral in the expression of the \textit{nir}-transform, but the second integral is ``discretised'', i.e. replaced by an infinite series.

$\exp_{\#}(X)$ (resp.  $\exp^{\#}(X)$) denotes the exponential function expanded as a power series of $X$ (resp.$X$ minus the leading term of $X$, which remains within the exponential). $\beta^{*}_{k} = \frac{B_{k+1}(\tau+\frac{1}{2})}{k+1}$ for the standard case, where $B_k$ is the $k^{\textrm{th}}$ Bernoulli polynomial and the function $f^{\Uparrow \beta^{*}}(n,\tau)$ and $f^{\uparrow \beta^{*}}(n,\tau)$ are defined as follows:
\begin{align*}
f^{\Uparrow \beta^{*}} &:= \beta({\partial}_{\tau})f(\frac{\tau}{n}) \quad \textrm{with}\quad {\partial}_{\tau}:=\frac{\partial}{\partial \tau} \\
&=:  f_{\kappa}\frac{n^{-\kappa} {\tau}^{k+1}}{k+1} + f^{\uparrow \beta^{*}}(n,\tau)
\end{align*}

\item The second, more analytical expression of \textit{nur}, consists in a succession of four (rather than nine) steps, three of which are elementary, but the fourth (the first, in fact) not at all, see section 5.2 [ES].

\item Then again \textit{nur} can be expressed as an infinite sum of \textit{nir}-transforms.
\[ nur(f)=\sum_{p \in \mathbb{Z}}(-1)^p nir(2 \pi i p + f) \]  

\end{enumerate}

\subsection{Singularities inferred from Taylor coefficient analysis.}
\bigskip
When dealing with  a SP series we are led to calculate the singularities using the asymptotics of its Taylor coefficients. Keeping the symmetry between the SP series $\phi$ and its Taylor coefficients $J$ intact, we look at them as resurgent functions in variables $z$ and $n$ and $\zeta$ and $\nu$ for different models that are explained in [ES]. For further details on this section [SS]
\\
The idea is to follow the correspondence between the triplets:
$$\{ \tilde{\phi}(z),\phi(\zeta),\phi(z) \} \leftrightarrow \{ \tilde{J}(n),J(\nu),J(n) \} $$
and that of the alien derivatives that are linked to them.

\bigskip
Retrieving the closest singularities:\\
Consider that our series $\hat{\phi}$ has one singularity at the boundary of the disk of convergence, which we denote by $\zeta_0$.
$$  \hat{\phi}(z) = \sum_{0 \leq n}(n+1)!J(n)z^{-n-1}(d{\nu}^t) $$ on applying the Borel transform gives $$ \hat{\phi}(\zeta) = \sum_{0\leq n} J(n) \zeta^{n} $$ 
To write down the closest singularity of $\hat{\phi}$ from the closest singularity of $\hat{J}$ , we write $J(n)$ as a cauchy integral on a circle $|\zeta| = |\zeta_0|-\epsilon $. 
$$ J(n) =\frac{1}{2 \pi i} \oint \hat{\phi}(\zeta)\zeta^{-n-1}d\zeta $$
By a contour deformation :
$$ =\frac{1}{2 \pi i} \int_{\Gamma} \hat{\phi}(\zeta)\zeta^{-n-1}d\zeta + o(\zeta^{-n}_0)$$
putting $\zeta := \zeta_0 \exp{\nu} = \exp^{\nu + \nu_{0}}$, we obtain 
$$ J(n) =\frac{1}{2\pi i} \int_{\Gamma_*} \hat{\phi}(\zeta_0\exp{\nu})\exp^{-n\nu}d \zeta + o(\exp^{-n\nu_{0}})$$

\noindent Using identities between $\hat{\phi}$ and $\check{\phi}$
$$ J(n) =\frac{1}{2\pi i} \int_{\Gamma_*} \check{\phi}_{\zeta_0}(\zeta_0 - \zeta_0\exp{\nu})\exp^{-n\nu}d \zeta + o(\exp^{-n\nu_{0}})$$

$$ J(n) =\frac{1}{2\pi i} \int_{\underline{\Gamma_*}} \check{\phi}_{\zeta_0}(\zeta_0 - \zeta_0\exp{-\nu})\exp^{n\nu}d \zeta + o(\exp^{-n\nu_{0}})$$
and in the case ${\stackrel{\diamond}{\phi}}_{\zeta_0}$ is integrable (that is to say in the case when the singularity at $\zeta_0$ is entirely recoverable from its minor i.e. free from Laurent terms), we arrive at:
$$\int^c_0 \hat{\phi}_{\zeta_0}(\zeta_0\exp{\nu}- \zeta_0)\exp^{-n\nu}d\zeta + o(\exp^{-n\nu_{0}})$$
In the above calculation we consider the contour $\Gamma$ which is the deformation of the circle $|\zeta|= |\zeta_0|-\epsilon$ to meet the larger circle $|\zeta|= |\zeta_0|+\epsilon$ avoiding the interval $[\zeta,\zeta + \epsilon]$.

 $\Gamma_*$ is the contour $\Gamma$ after undergoing $\zeta =\zeta_0\exp{\nu}$.
The other distant singularities could be obtained by extending this procedure. An important remark here is that if $\hat{J}$ is endlessly continuable, its resurgence pattern can be very precisely obtained from that of the $\hat{\phi}$.

%% file: ode_2.tex
\section{Differential aspects}

The inner generators for some driving functions (for instance all polynomial $f$ and all monomial $F$) satisfy ordinary differential equation of linear homogeneous types with polynomial coefficients of usually high order. These equations can be particularily interesting as they give a rather explicit explanation of the inner generators which occur in our problem in an algebraic manner. They enable a detailed understanding of their behaviour at $\infty$ in the $\nu$ plane or at the origin in the Borel $\zeta$-plane. 

What is a bit more appealing than the rest is the occurence of the phenomenon of \textit{rigid} \,\textit{resurgence} which involves discrete Stokes constants in the asymptotic expansions. The reason it might be called rigid is that it exhibits a certain stability to the changing parameters of the problem: contrary to what happens with most singular ODEs, here the resurgence coefficients are discrete; they don't vary continously with the real or complex parameters inside the ODE. More detail can be found in [ES]

In the developpement of the differential properties of the inner generators, we would talk about the so called ``variable'' and the ``covariant'' differential equations. The ``variable'' equation applies to the whole \textit{nir}-transform, whereas the ``covariant'' equation applies to the non-entire part (usually semi-entire part) of the \textit{nir}-transform, which alone has intrinsic geometric meaning: it describes the singularity present at this point. As a consequence, the ``variable'' equation is specific to inner generator, whereas the ``covariant'' equation applies equally to all inner generators under a shift of base point. In what follows, we would follow closely the sketch drawn for the treatment of this problem in chapter 6 of [ES]. The notations stay the same and to illustrate the problem and to be able to handle this section efficiently, we would like to write down a reminder of the basics.
 

Let's start with writing down the four types of shift operators $ \beta ({\partial}_{\tau})$ to settle notation:
\begin{eqnarray*}
\textrm{trivial choice} \,\,\beta ({\tau}) &:=& {\tau}^{-1} \\
\textrm{standard choice} \,\,\beta ({\tau}) &:=& ({\exp}^{{\tau}/2} -  {\exp}^{{-\tau}/2})^{-1} \\
\textrm{odd choice} \,\,\beta ({\tau})&:=& {\tau}^{-1} + \sum_{s\geq 0} \beta_{2s+1}\tau^{2s+1}\\
\textrm{general choice} \,\,\beta ({\tau})&:=& \sum_{s\geq 0} \beta_{s}\tau_{s}
\end{eqnarray*}
Consider the driving function $f$. It can be subjected to the nir-transform provided $f(0)=0$ and more particularily $f^{\prime}(0) \neq 0$.
\begin{eqnarray*}
 f(x)&:=& \sum_{1\leq s \leq r} f_s x^s \\
\phi(n,\tau)&:=& \beta(\tau)f(\frac{\tau}{n}) = \frac{1}{2}\frac{{\tau}^2}{n}f_1 + \dots \\
\phi(n,\tau)&:=& \phi^{+}(n,\tau) + \phi^{-}(n,\tau) \\
\phi^{\pm}(n,\pm \tau)&:=& \pm \phi^{\pm}(n,\pm \tau)\\
k(n) &:=& [ \int^{\infty}_0 \exp ^{\#}(\phi(n,\tau)) d\tau ]_{\textrm{singular}}\\
&:=& \int^{\infty}_0 \exp ^{\#}(\phi(n,\tau))cosh_{\#}(\phi^-(n,\tau)) d\tau \,\,\,(\textrm{if} f_1 \neq 0)\\
k(n) &:=& [\frac{1}{2\pi\iota} \int^{c+\iota\infty}_{c-\iota\infty} k(n) \exp^{ \nu n} dn ]_{\textrm{formal}} = h^\prime(\nu)\\
k(n) &:=& [\frac{1}{2\pi\iota} \int^{c+\iota\infty}_{c-\iota\infty} k(n) \exp^{\nu n} \frac{dn}{n} ]_{\textrm{formal}} = h(\nu)
\end{eqnarray*}
\\
\subsection{Variable and covariant linear homogeneous polynomial ODEs }
The general form of such equations :
$$ \textrm{variable ODE}: P_v (n, -\partial_n) k^{total}(n) = 0 \Leftrightarrow P_v (\partial_{\nu},\nu)k^{total}(\nu)=0 $$ 
$$ \textrm{covariant ODE}: P_c (n, -\partial_n) k^{singular}(n) = 0 \Leftrightarrow P_c (\partial_{\nu},\nu)k^{singular}(\nu)=0 $$
where:
$$ P_v (n,\nu) = \sum_{0 \leq p \leq d} \sum_{0 \leq q \delta} dv_{p,q}n^p \nu^\psi_k (n,\tau)q $$
$$ P_c (n,\nu) = \sum_{0 \leq p \leq d} \sum_{0 \leq q \delta} dc_{p,q}n^p \nu^q $$

These polynomials are of degree $d$ and $\delta$ in the variables $n$, $\nu$. Here $n$ and $\nu$ should be treated as non-commutative variables subject to the following relation $[n,\nu]=1$.
\newline
The relation of covariance is as follows: $$ P^{^\epsilon f}_{c}(n,\nu-\eta)\equiv
P^f_c (n,\nu) \,\,\,\,\,\textrm{for all}\,\, \epsilon $$
with: $ ^\epsilon f(x) = f(x+\epsilon)$ and $ \eta:= \int^{\epsilon}_0 f(x) dx $

Here we give a little note about the existence and the calculation of these ODEs. The details can be found [ES]. The polynomials $ \varphi_s $ and $\psi_s $ in $ n^{-1} $ and $\tau$ for any $ s \in \mathbb{N} $ can be written down using the following identities:
$$ \partial^s_n k^{total}(n) = \int^{\infty}_0 \varphi_s (n,\tau) \exp^{\#}(\varphi(n,\tau)) d\tau $$
with $$ \varphi_s (n,\tau) \in \mathbb{C}[\partial_n\varphi,\partial^2_n\varphi, \dots,\partial^s_n\varphi] \in \mathbb{C}[{n^{'}-1}, \tau]$$
$$ \int^{\infty}_0 d^s_{\tau}(\tau^k \exp^{\#}(\varphi(n,\tau))) = \int^{\infty}_0 \psi_s(n,\tau) \exp^{\#}(\varphi(n,\tau)) d\tau = 0 $$
with $$ \psi_s(n,\tau) = \tau^s \partial_{\tau}\varphi(n,\tau)+ s\tau^{s-1} \in \mathbb{C}[{n^{-1}, \tau}] $$
\newline
For larger degrees the polynomials $\varphi_s$and $\psi_s$ become linearly independant on $\mathbb{C}[n^{-1}]$. 
$$ \sum_{0 \leq s \leq \delta} A_s(n)\varphi(n, \tau) + \sum_{0 \leq s \leq \delta^{'}} B_s(n)\varphi(n, \tau) = 0 \,\,\textrm{for}\,\, A(n), B(n) \in \mathbb{C}[n] $$
and to each such relation there corresponds a linear ODE for $k^{total}$ :
$$ ( \sum_{0 \leq s \leq \delta} A_s(n)\partial^s_n) k^{total}(n) = 0 $$

The expressions given here work for $f(0)=0$ also. In this case when we apply the \textit{nir} transform to the corressponding driving function, we obtain powers that are integer as well as fractional and hence the use of the term $k^{total}$ which implies that one considers all powers for an overall analysis of the inner generators. In this sction we talk about the polynomial $f$ inputs for which $k$ and $k^{total}$ in both $x$ and $\nu$ plane, satisfy ODEs. 

The ones satisfied by $k^{total}$ are called variable as they strictly depend on a proper base point. 
Whereas the ones satisfied by $k$ are called covariant: for a change of a proper base point $x_i$ to $x_j$, the covariant differential equation satisfied in the $\nu$ plane undergoes a translation in the $\nu$ plane of $\nu = \int^{x_i}_{x_j}f(x) dx $. Here $x_0$ is proper base point means that $f(x_0)=0$.  Another important thing about the covariant ODEs is that there exists a unique extension of the covariant ODE to a non proper base point $x_l$. This evidently does not coincide with the variable ODE.


\subsection{Method for calculating the differential equation}

A rather schematic version of the entire problem. For details for explicit calculations with the technical aspects, the reader is invited to refer to [ES]. There are detailed reports on schemes on the existence and calculation of variable ODEs and the covariant ODEs ($f(0)=0$ and $ f(0) \neq 0$).

The integrals for the inner generators: $h$, $h^{\prime}$ and $hh$ are as follows:
\begin{eqnarray*} 
h(\nu) &=& \frac{1}{2 \pi i} \int^{+i \infty}_{ -i \infty} e^{n \nu} \frac{dn}{n} \int^{\infty}_0 e^{\phi (n,\tau)} d\tau \,\,\,\,\,\textrm{(singular germ)}\\
h^{\prime}(\nu) &=& \frac{1}{2 \pi i} \int^{+i \infty}_{ -i \infty} e^{n \nu} dn \int^{\infty}_0 e^{\phi (n,\tau)} d\tau \,\,\,\,\,\textrm{(singular germ)}\\
hh(n) &=& \int^{\infty}_0 e^{\phi (n,\tau)} d\tau \,\,\,\,\,\,\,\,\,\,\,\,\,\,\textrm{(divergent power series)}\\
\end{eqnarray*}
for $ \phi (n,\tau):= - \beta({\partial_{\tau}}).f(\frac{\tau}{n})$.
\\
The four choices for $\beta$, the shift operator, are as has been mentioned before in the article, namely trivial,  standard, odd and general.
\\
Now one can write down the series for $h(\nu)$, $h^{\prime}(\nu)$ and $hh(n)$:

\begin{align*}  h(\nu) & = \sum_{m \in {\mathbb{N}}_k} h_m {\nu}^m \\
 h^{\prime} (\nu) & = \sum_{m \in {\mathbb{N}}_k} h_m  m {\nu}^{m-1} \\
 hh(n) & = \sum_{m \in {\mathbb{N}}_k}  h_m  m! n^{-m} \end{align*}

with ${\mathbb{N}}_k := -\frac{\kappa}{1+\kappa} + \frac{1}{1+\kappa} \mathbb{N}$ \\


\subsubsection{The Differential Equations:}

The inner generators $h$, $h^{\prime}$ and $hh$ verify the differential equations $D$, $D^{\prime}$ and $DD$ respectively. (These differential equations are of linear and homogeneous type and fall into both, the variable and covariant category).

$$ D = \sum_{\substack{ 
0 \leq d \leq d^*, \\
0 \leq \delta \leq \delta^*}} D_{d,\delta} \nu^d \partial^{\delta}_{\nu} $$ 

$$ D^{\prime} = \sum_{\substack{ 
0 \leq d \leq d^*, \\
0 \leq \delta \leq \delta^*}} D^{\prime}_{d,\delta} \nu^d \partial^{\delta}_{\nu} $$ 

$$ DD = \sum_{\substack{ 
                   0 \leq d \leq d^*, \\
                   0 \leq \delta \leq \delta^*}} D_{d,\tau} n^{-\delta} \partial^{d}_{n} \,= \,  \sum_{\substack{ 
0 \leq d \leq d^*, \\
0 \leq \delta \leq \delta^*}} D_{d,\delta} ({{-\partial}_n})^{d} n^{-\delta}$$  

\noindent The last equality uses the equivalence $({\partial}_{\nu}, -\nu) \leftrightarrow (n^{-1},{\partial}_{n})$ \\ 


\noindent One observes the following:
$$ {({\partial}_n)^k} hh(n) \equiv \int^{\infty}_0 \phi_k (n,\tau) e^{\phi_k (n,\tau)} d\tau  \,\,\,\,(\forall k \in \mathbb{N}) $$
$$ 0 \equiv \int^{\infty}_0 \partial_{\tau}[{\tau}^k e^{\phi_k (n,\tau)}] d\tau = \int^{\infty}_0 \psi_k (n,\tau) e^{\phi_k (n,\tau)} d\tau  \,\,(\forall k \in \mathbb{N}^*) $$

\noindent where $\psi_k (n,\tau) := {\tau}^k \partial_{\tau} \phi_k (n,\tau) + k {\tau}^{k-1} $

Then for large enough values of the $d^*$ and $\delta^*$, the elementary functions $ \phi_k (n,\tau)$ for $(k= 0,1,\dots,d^*)$ and $\psi_k (n,\tau)$ for $(k= 0,1,\dots,{\delta}^*)$, 	are no longer linearly independant. We use this fact to eliminate the $\psi_k (n,\tau)$ to get an ODE that is satisfied by $hh(n)$.

\textbf{Check}: To verify that a certain test function satisfies an ODE, the way we go about it is to explicitely calculate the nir transform for the function and solve for a differential system with a convenient number of successive derivatives. This is the  method that we have used extensively in the course of this article to prove the non existence of elementary (i.e. low-order and low degree) differential equations (even algebraic) for certain class of driving functions $f$ (or $F$).
\bigskip
\subsubsection{Dimensions}
Dimensions of the spaces of variable ODEs. 

For $r:= \textrm{def}(f)$ and for each pair $(x,y)$ with
$$ x \in \{ v,c \} = { variable, covariant} $$
$$ y \in \{ t,s,o,g\} = \{ trivial, standard, odd, general \} $$ 
It turns out that in each case, the dimension of the space of the ODEs is of the form:
$$ dim_{x,y}(r,d,\delta) \equiv (d- A_{x,y}(r))(\delta- B_{x,y}(r))-C_{x,y}(r) $$
for $\delta$ and $d$, the differential order of the ODEs in the $n$- variable (and d for the $\nu$ -variable). The extremal pairs $(\underline{d},\ - \delta)$ and $( -d, \underline{\delta})$ for 
$$ \underline{d} = 1 + A_{x,y}(r)$$
$$ \-{d} = 1 + A_{x,y}(r) +  C_{x,y}(r)$$
$$ \underline{\delta} = 1 + B_{x,y}(r)$$
$$ \-{\delta} = 1 + B_{x,y}(r) + C_{x,y}(r)$$
($\underline{d}$ and $\underline{\delta}$ is minimal, $\-{d}$ and $\-{\delta}$ is cominimal), the corresponding dimension is always exactly one.


The details for the dimensions of the variable and covariant ODEs in the various cases (trivial, standard, odd and general) are specified in the article [ES].
\\

\subsubsection{The equivalence: $({\partial}_{\nu}, -\nu) \Longleftrightarrow (n^{-1},{\partial}_{n})$}

The differential operator for both the variable ODEs and the covariant ODEs, can be expressed as polynomials $P(n,\nu)$ with degree $ (d, \delta ) $. They are written down as functions in $(n,\nu)$. 

$(n,\nu)$ are a non commutating pair of variables with $[n,\nu]=1$. Hence there are two equivalent parametrisations possible (consistent with the relation of non-commutation)
$$ (n,\nu)\,\,\Longrightarrow\,\, (n,-\partial_{n})\,\,\,\,\,\,(\star)$$
$$ (n,\nu)\,\,\Longrightarrow\,\, (\partial_{\nu},\nu)\,\,\,\,\,\,(\star\star)$$

The corresponding ODE equivalence can be written as:
$$ P(n,-\partial_n) k(n) = 0 \,\,\Longleftrightarrow \,\,  P(\partial_{\nu}, \nu) \hat{k}(n) = P(\partial_{\nu}, \nu)\partial_{\nu}h(\nu) = 0 $$

%

\subsubsection{Compressing the covariant ODEs}
This is essentially the use of the covariance relation to make the ODEs easier to handle.
The three step procedure for this:

 \begin{enumerate}
\item Apply the ODE algorithm to the centered function $f(x)= \sum^{r-2}_{i=0} f_i x^i + f_r x^r $.

\item The coefficients $\{f_0,f_1,\dots,f_{r-2},f_r\}$ are replaced by the shift-invariants  $\{\textbf{f}_0,\textbf{f}_1,\dots,\textbf{f}_{r-2},\textbf{f}_r\}$.

\item Then the $\beta$-coefficients are replaced by the centered $\bm{\beta}$-coefficients.
\end{enumerate}

\subsubsection{Basic polynomials $f(x)$ and $p(\nu)$}
$$ f(x) = f_0+f_1 x + \dots + f_r x^r = (x-x_1)\dots(x-x_r)f_r $$
$$ p(\nu) = p_0+p_1 \nu + \dots + p_r {\nu}^r = (\nu-{\nu}_1)\dots(\nu-{\nu}_r)p_r $$
with
$$ \nu_i = f^{*}(x_i) = \int^{x_i}_0 f(x)dx = \sum_{0\leq s \leq r} f_s \frac{x^{s+1}_i}{s+1} $$

\subsubsection{Basic symmetric functions}
$$ \{f_s\} \rightarrow \{x^*_s\} \rightarrow \{ x^{**}_s\} \rightarrow \{{\nu}^{**}_s\} \rightarrow \{{\nu}^*_s\} \rightarrow \{p_s\} $$
\begin{itemize}
 \item $$ \sum_{1 \leq s \leq \infty } \frac{1}{s} \frac{x^{**}_s}{x^s} = -{\log}(1+\sum_{1 \leq s \leq r} {(-1)}^r \frac{x^{*}_s}{x^s}) $$

 \item $$ {\nu}^{**}_s \equiv \sum_{s \leq t \leq (r+1)s} f^*_{s,t}x^{**}_t \,\,\textrm{with}\, \,\sum_{s \leq t \leq (r+1)s} f^*_{s,t}x^t := (f^*(x))^s $$

  \item $$ 1+(\sum_{1 \leq s \leq r} {(-1)}^r \frac{{\nu}^{*}_s}{{\nu}^s}) \equiv \exp(-\sum_{1 \leq s \leq \infty} \frac{{\nu}^{**}_s}{{\nu}^s}) $$
\end{itemize}

\subsection{Covariant ODEs and leading covariant polynomials}
\subsubsection{Covariant ODEs}

We consider the inner generators associated with the driving function. In the Borel plane $\nu$, they are singular germs $h(\nu)$. In the multiplicative plane $n$, they are divergent power series $hh(n)$ in decreasing powers of $n$. Any linear homogeneous ODE for $h(\nu)$ translates into a linear homogeneous ODE for $hh(n)$, and vice-versa, under the change $(\partial_{\tau},\nu) \leftrightarrow (n,-\partial_{n})$, in keeping with $(\star)$, $(\star \star)$ 
$$ CD. h(\nu) = 0 \,\,\,\,\,\, CDD. hh(n) = 0 $$
The shape of these equations is linear homogeneous. The covariance of the differential equation is given as before.

Covariant differential equation has 
\begin{itemize}
 \item invariant coefficients 
$$ A: f \mapsto A(f) \,\, \textrm{with}\,\, A(f_{\epsilon}) \equiv A(f) \,\,\, \forall \epsilon $$
 \item covariant polynomials
$$  P: f \mapsto P(f, \nu) \,\, \textrm{with}\,\, P(f_{\epsilon},\nu) \equiv  P(f,\nu + \eta) \,\,\, \forall \epsilon $$
with $\eta := \int^{\epsilon}_{0}f(x) dx $ and $ f_{\epsilon}(x):= f(x+ \epsilon)$
\end{itemize}

\subsubsection{Leading covariant polynomials $\Lambda(\nu)$}

For :
$$f(x):= \sum_{0\leq i \leq r} f_i x^i = f_r \prod^r_{i=1} (x-x_i) \mapsto \Lambda(\nu) = f_r \prod^r_{i=1} (\nu-\nu_{i})$$
here $\nu_{i} :=  \int^{x_i}_{0} f(x) dx$
\begin{itemize}
 \item By the covariance of $P(f, \nu) := \Lambda(\nu)$
\item If we write  $$ CD^{f} (\nu, {\partial}_{\nu}) = \sum_{0 \leq \delta \leq {\delta}^*} CD^f_{\delta}(\nu){\partial}^{\delta}_{\nu} $$
Then the leading term is $ CD^f_{{\delta}^*}(\nu) $ (term in front of the leading derivative ${\partial}^{{\delta}^*}_{\nu} $). 
\\
is a covariant polynomial 
\\
is a multiple of $\Lambda(\nu)$.

\end{itemize}

The covariant polynomial $\Lambda(\nu)$ is important as it is responsible for the singularities of $h(\nu)$.

\subsubsection{Characterisation of the invariance and covariance}
For an $r$ (degree of $f(x)$) fixed, one sets:

 $$ {\partial}_{\epsilon} f_{0} := f_{1}; {\partial}_{\epsilon} f_{1} := 2.f_{2};
{\partial}_{\epsilon} f_{2} := 3.f_{3};\dots; {\partial}_{\epsilon} f_{i} :=(i+1)f_{i+1} ;\dots; {\partial}_{\epsilon} f_{r} := 0 $$
$${\partial}_{\epsilon} \nu := -f_{0}$$

Then;
$$ A: f \longrightarrow A(f)= A(f_0,f_1,\dots,f_r) \,\, \textrm{\,is\,invariant\, iff}\,\,{{\partial}_{\epsilon}}A(f_0,f_1,\dots,f_r) \equiv 0 $$

$$ P: f \longrightarrow P(f,\nu)= O(f_0,f_1,\dots,f_r;\nu) \,\, \textrm{\,is\,covariant\, iff}\,\,{{\partial}_{\epsilon}}P(f_0,f_1,\dots,f_r;\nu) \equiv 0 $$

\subsubsection{Invariants}
Here are some invariants for various values of $r$ for illustration.
\noindent
\\
$\,r=1: \,\,\, f_1$
\\
$\,r=2: \,\,\, f_2;\,\,\,f_1f_1-4.f_0f_2$
\\ 
$\,r=3: \,\,\, f_3;\,\,\,f_2f_2-3.f_1f_3 \,\,\, f_1f_1f_2f_2-4f_1f_1f_1-4f_0f_2f_2f_2+18f_0f_1f_2f_3-27f_0f_0f_3f_3 $
\\
\subsubsection{Canonical covariant polynomials of degree 1}
$\,r=1:\,\,\Pi_1(\nu)\,=\,\,f_1\nu + \frac{1}{2}f_0f_0$
\\
$\,r=2:\,\,\Pi_2(\nu)\,=\,\,f_2f_2\nu + \frac{1}{2}f_0f_1f_2 -\frac{1}{20}f_1f_1f_1 $
\\
$\,r=3:\,\,\Pi_3(\nu)\,=\,\,f_3\nu + \frac{1}{3}f_0f_2 -\frac{1}{20}f_1f_1 $
\\
$\,r=4:\,\,\Pi_4(\nu)\,=\,\,f_4f_4\nu + \frac{1}{4}f_0f_3f_4 -\frac{1}{60}f_1f_2f_4 -\frac{1}{40}f_1f_3f_3+ \frac{1}{120}f_2f_2f_3$
\\
$\,r=5:\,\,\Pi_5(\nu)\,=\,\,f_F\nu + \frac{1}{5}f_0f_4 -\frac{1}{20}f_1f_3 + \frac{1}{60}f_2f_2 $
\\
Here we remark that there is an alternate double homogeneity which depends on the parity of $r$. (For odd values the sum terms is $[r-1,2]$ and for even values $[2r-1,3]$).


%% file: ode_7.tex
\subsection{The global resurgence picture for the polynomial inputs $f$}
Using the covariant ODEs one could describe the exact singular behaviour of $ k( \nu) = h(\nu)$ at infinity in the $\nu$-plane, and so rather evidently the singularities at zero in the $\zeta $-plane. The solutions of the covariant ODE at infinity , in the $\nu$-plane, are linear combinations of elementary exponential factors along with the negative powers of $\nu$. They are always divergent, resurgent and resummable.
\\
For radial inputs like $ f(x) = f_r x^r $: the power series is of the form:
$$ \sum_{r+1 \leq k} c_s(\omega)\nu^{-\frac{s}{r+1}} exp(\omega \nu^{\frac{r}{r+1}} )$$

For a general $ f $, it is of the form:

$$ \sum_{r+1 \leq k} c_s(\omega)\nu^{-\frac{s}{r+1}} exp(\omega \nu^{\frac{r}{r+1}} + \sum^{r-2}_{s=1}{\omega}_s \nu^{\frac{s}{r+1}} ) $$

with ``frequencies'' $\omega_{s}$ that depend algebraically on the co-efficients of the ODEs.

\section{Rational inputs $F$ : the resurgence pattern and connection matrices.}

\subsection{Matrice entries}

Let us display the matrices of transformation which express the linear transformations undergone by the Laplace transforms (of direction $\theta$, of all inner generators) when the integration axis $\theta$ crosses a singular axis $\theta_{0}$, i.e. a singularity carrying axis.
The various steps that have been taken into consideration :
after the figures for the cases $ p= 3, \dots , 9 $ and the transformation vectors for the various cases , we would like to write down explicitly the transformation matrices with the other intermediate matrices that are part of the calculations. We specifically consider the case for half , full , three-fourth and two turns in the plane of the singularities and we would constrict ourselves to the case of $ \epsilon = -1 $. For each case, for at least the lower values of $p$ (till $p = 9$ ) we give the transpose of the matrix $M$ as this is more directly related to our problem, and for higher $p$ values we give just the characteristic polynomial of the matrix.

$$\mathcal{R}_{{\theta}_{\frac{1}{2}}}: = \left(
\begin{array}{c}
x_1\\
x_2  
\end{array}
\right) \rightarrow \left(
\begin{array}{c}
x_1 \\
x_2 + 2x_1{\epsilon}^\frac{1}{2}
\end{array}
\right) \\ \mathcal{R}_{{\theta}_{1}}: = \left(
\begin{array}{c}
x_1\\
x_2 
\end{array}
\right) \rightarrow \left(
\begin{array}{c}
x_1 + 2x_2{\epsilon}^\frac{1}{2} \\
x_2
\end{array}
\right) $$

\bigskip

\begin{tikzpicture}
  \node at (0,0) (O) {$\bullet$};
\draw [very thin](0,0) circle(5cm);
\draw [dotted](0:5)  \foreach \x in {120,240} {
                -- (\x:5) 
            }-- cycle (-90:5) node[below=1cm] {$n=3$} ;
\draw [->>,thick](180:6)--(0:6) (0:5) node[below=2mm] {$x_0 = x_3$};
\draw [->>,thick](300:6) -- (120:6) (120:5) node[left=2mm]{$x_1$};
\draw [->>,thick](240+180:6) -- (240:6) (240:5) node[left=2mm]{$x_2$};

\end{tikzpicture}
\\

$$\mathcal{R}_{{\theta}_{\frac{1}{2}}}: = \left(
\begin{array}{c}
x_1\\
x_2 \\
x_3 
\end{array}
\right) \rightarrow \left(
\begin{array}{c}
x_1 \\
x_2 + 3x_3{\epsilon}^\frac{1}{3}  \\
x_3 
\end{array}
\right) \\ \mathcal{R}_{{\theta}_{1}}: = \left(
\begin{array}{c}
x_1\\
x_2 \\
x_3 
\end{array}
\right) \rightarrow \left(
\begin{array}{c}
x_1 \\
x_2 + 3x_1{\epsilon}^\frac{2}{3}\\
x_3 
\end{array}
\right) $$
$$\mathcal{M}_{\frac{1}{2}}: = \left(
\begin{array}{ccc}
 1 & 3 {\epsilon}^{\frac{2}{3}} & 3 {\epsilon}^{\frac{1}{3}} \\
 0 & 1                          & 0                          \\
 0 & 3 {\epsilon}^{\frac{1}{3}} & 1
\end{array}
\right) \mathcal{M}_{1}: = \left(
\begin{array}{ccc}
 10                      & 3 {\epsilon}^{\frac{2}{3}} & -6 {\epsilon}^{\frac{1}{3}}\\
 -6 {\epsilon}^{\frac{1}{3}} & 1                      &  3 {\epsilon}^{\frac{2}{3}}\\
 -15 {\epsilon}^{\frac{2}{3}}& 3 {\epsilon}^{\frac{1}{3}} & -8
\end{array}
\right) $$


The characteristic polynomials for ${\epsilon}$ for the intermediate matrices are the following when we turn around once:
$$ (t-1)^3 , (t-1)^3 ,(t-1)^3 ,(t-1)(t^2 + 7t +1), t^3 + 15 t^2 +12t -1 , (t-1)^3 $$

\begin{tikzpicture}
  \node at (0,0) (O) {$\bullet$};
\draw [very thin](0,0) circle(5cm);
\draw [dotted](0:5)  \foreach \x in {90,180,...,360} {
                -- (\x:5) 
            }-- cycle (-90:5) node[below=1cm] {$n=4$};
\draw [->>,red,thick](180:6)--(0:6) (0:5) node[below=2mm] {$x_0 = x_4$};
\draw [->>,thick](270:6) -- (90:6) (90:5) node[above=2mm] {$x_1$}; 
\draw [->>,thick](360:6) -- (180:6) (180:5) node[below=2mm] {$x_2$};
\draw [->>,thick](270+180:6) -- (270:6) (270:5) node[below=2mm] {$x_3$};

\end{tikzpicture}

$$\mathcal{R}_{{\theta}_{\frac{1}{2}}}:= \left(
\begin{array}{c}
x_1\\
x_2 \\
x_3 \\
x_4
\end{array}
\right) \rightarrow \left(
\begin{array}{c}
x_1 \\
x_2 + 4x_1{\epsilon}^\frac{3}{4} \\
x_3 + 4x_4{\epsilon}^\frac{1}{4} \\
x_4
\end{array}
\right) \\ \mathcal{R}_{{\theta}_{1}}:= \left(
\begin{array}{c}
x_1\\
x_2 \\
x_3 \\
x_4
\end{array}
\right) \rightarrow \left(
\begin{array}{c}
x_1 \\
x_2 \\
x_3 + 6x_1{\epsilon}^\frac{1}{2} \\
x_4
\end{array}
\right)$$
$$\mathcal{M}_{\frac{1}{2}}: = \left(
\begin{array}{cccc}
 1&4{\epsilon}^{\frac{3}{4}}&-10\iota&-20{\epsilon}^{\frac{1}{4}}\\
0&1&4{\epsilon}^{\frac{3}{4}}&6\iota \\
0&0&1&0\\
0&0&{\epsilon}^{\frac{1}{4}}&1
\end{array}
\right) $$

$$ \mathcal{M}_{1}: = \left(
\begin{array}{cccc}
45& 4(-20{\epsilon}^{\frac{1}{4}} - 40 \iota {\epsilon}^{\frac{3}{4}}){\epsilon}^{\frac{3}{4}}& 4{\epsilon}^{\frac{3}{4}} &-40 \iota {\epsilon}^{\frac{1}{4}}\\ 
6 \iota (-20{\epsilon}^{\frac{1}{4}} - 40 \iota {\epsilon}^{\frac{3}{4}}) & -10 \iota & (-20{\epsilon}^{\frac{1}{4}} - 40 \iota {\epsilon}^{\frac{3}{4}})& 36
\\
-10 \iota & -20{\epsilon}^{\frac{1}{4}} & 1& 4{\epsilon}^{\frac{3}{4}}\\
-36{\epsilon}^{\frac{3}{4}}&-74 \iota & 4{\epsilon}^{\frac{1}{4}}&-15
\end{array}
\right)
$$

$$\mathcal{M}_{\frac{3}{2}}: = \left(
\begin{array}{cccc}
45 + 4h{\epsilon}^{\frac{3}{4}} & r & -10\iota +4h{\epsilon}^{\frac{1}{4}} + 6\iota(45 + 4h{\epsilon}^{\frac{3}{4}} +4 r{\epsilon}^{\frac{3}{4}}) &  h + 4(45+4h({\epsilon}^{\frac{1}{4}})+6\iota r)\\
h & 45+ 4h {\epsilon}^{\frac{3}{4}}& r & -10\iota + 4{\epsilon}^{\frac{1}{4}} + 6\iota(45 + 4h{\epsilon}^{\frac{3}{4}})\\
-10\iota & h & 45+ 4h {\epsilon}^{\frac{3}{4}}& 4{\epsilon}^{\frac{3}{4}}-40 \iota {\epsilon}^{\frac{1}{4}} + 6 \iota h\\
-36{\epsilon}^{\frac{3}{4}}& 70\iota& -120{\epsilon}^{\frac{1}{4}}& -291
\end{array}
\right)
$$
for $r$ = $ 4(-1)^{\frac{3}{4}} -40 \iota (-1)^{\frac{1}{4}} + 6\iota (-20(-1)^{\frac{1}{4}} - 40 \iota (-1)^{\frac{3}{4}}) +4(45+4(-20(-1)^{\frac{1}{4}} - 40 \iota (-1)^{\frac{3}{4}})((-1)^{\frac{3}{4}}))(-1)^{\frac{3}{4}}$ and $h=-20{\epsilon}^{\frac{1}{4}} - 40 \iota {\epsilon}^{\frac{3}{4}}$

 


The characteristic  polynomials for the intermediate matrices (in the case of even $p$ we would consider two turns instead of a single one).
$$ (t - 1)^4 , (t - 1)^4 , (t - 1)^4 , (t - 1)^4 ,$$
$$ (t^2  + 14 t + 1)^2 ,1 - 32 t - 258 t^2  - 32 t^3  + t^4 ,$$
$$ 1 + 64 t - 258 t^2  + 64 t^3  + t^4 , (t + 1)^4 ,$$
$$ (t^2  + 106 t + 1) (t^2  - 6 t + 1), 1 + 40 t - 370 t^2  + 40 t^3  + t^4 ,$$
$$ 1 - 24 t + 590 t^2  - 24 t^3  + t^4 , 1 + 396 t + 3494 t^2  + 396 t^3  + t^4 ,$$
$$ 1 - 308 t + 5286 t^2  - 308 t^3  + t^4 , 1 + 196 t - 4746 t^2  + 196 t^3  + t^4 ,$$
$$ 1 - 508 t - 2058 t^2  - 508 t^3  + t^4 , (t - 1)^4 $$

\begin{tikzpicture}
  \node at (0,0) (O) {$\bullet$};
\draw [very thin](0,0) circle(5cm);
\draw [dotted](0:5) 
\foreach \x in {72,144,...,359} {
            -- (\x:5)
        } -- cycle  (-90:5) node[below=1cm] {$n=5$} ;

\draw [->>,thick](180:6) -- (0:6) (0:5) node[below=2mm] {$x_0 = x_5$}; 
\draw [->>,thick](-108:6) -- (72:6) (72:5) node[above=2mm] {$x_1$}; 
\draw [->>,thick](144+180:6) -- (144:6) (144:5) node[left=2mm] {$x_2$};
\draw [->>,thick](-144+180:6) -- (-144:6) (-144:5) node[below=2mm] {$x_3$};
\draw [->>,thick](-72+180:6) -- (-72:6) (-72:5) node[below=2mm] {$x_4$};
[>=stealth]
  \draw[<<->>,thick]
    (10mm,0) arc (0:31:10mm); node[above=1cm,fill= blue]{$ \theta $};

[>=stealth]
\draw[<<->>,thick,blue!80]
    arc (31:72:10mm);


\end{tikzpicture}

$$\mathcal{R}_{{\theta}_{\frac{1}{2}}}: = \left(
\begin{array}{c}
x_1\\
x_2 \\
x_3 \\
x_4 \\
x_5
\end{array}
\right) \rightarrow \left(
\begin{array}{c}
x_1 \\
x_2 + 5x_1{\epsilon}^\frac{4}{5} \\
x_3 + 10x_5{\epsilon}^\frac{2}{5} \\
x_4 \\
x_5

\end{array}
\right) \\ \mathcal{R}_{{\theta}_{1}}: = \left(
\begin{array}{c}
x_1\\
x_2 \\
x_3 \\
x_4 \\
x_5
\end{array}
\right) \rightarrow \left(
\begin{array}{c}
x_1 \\
x_2 + 5x_1{\epsilon}^\frac{4}{5} \\
x_3 + 10x_5{\epsilon}^\frac{2}{5} \\
x_4 \\
x_5

\end{array}
\right)$$
$$\mathcal{M}_{\frac{1}{2}}: = \left(
\begin{array}{ccccc}
 1 & 5{\epsilon}^{\frac{4}{5}} & -15{\epsilon}^{\frac{3}{5}}  & 35{\epsilon}^{\frac{2}{5}}  & 45{\epsilon}^{\frac{1}{5}}\\
 0 & 1                   &  5{\epsilon}^{\frac{4}{5}}   & -15{\epsilon}^{\frac{3}{5}} & 10{\epsilon}^{\frac{2}{5}}\\
 0 & 0                   &  1                     &  5{\epsilon}^{\frac{4}{5}}  & 0                   \\
 0 & 0                   &  0                     &  1                    & 0                   \\
 0 & 0                   &  10{\epsilon}^{\frac{2}{5}}  & -45{\epsilon}^{\frac{1}{5}} & 1

\end{array}
\right)$$
$$ \mathcal{M}_{1}: = \left(
\begin{array}{ccccc}
 126                   & -420{\epsilon}^{\frac{4}{5}} & 160{\epsilon}^{\frac{3}{5}}    & 35{\epsilon}^{\frac{2}{5}}  & -70{\epsilon}^{\frac{1}{5}}\\
 -70{\epsilon}^{\frac{1}{5}} & -224                   &  -70{\epsilon}^{\frac{4}{5}}   & -15{\epsilon}^{\frac{3}{5}} & 35{\epsilon}^{\frac{2}{5}} \\
 35{\epsilon}^{\frac{2}{5}}  & 105{\epsilon}^{\frac{1}{5}}  &  -24                     &  5{\epsilon}^{\frac{4}{5}}  & -15{\epsilon}^{\frac{3}{5}} \\
 -15{\epsilon}^{\frac{3}{5}} & -40{\epsilon}^{\frac{2}{5}}  &  5{\epsilon}^{\frac{1}{5}}     &  1                    & 5{\epsilon}^{\frac{4}{5}}   \\
 280{\epsilon}^{\frac{4}{5}} & 860{\epsilon}^{\frac{3}{5}}  & -215{\epsilon}^{\frac{2}{5}}   & -45{\epsilon}^{\frac{1}{5}} & 126
\end{array}
\right)
$$

The characteristic polynomials for all the intermediate matrices once again for ${\epsilon} = -1$ are as follows:
$$ (t-1)^5 , (t-1)^5 , (t-1)^5 , (t-1)^5, (t-1)^5,$$
$$ (t - 1) (t^2  + 23 t + 1) (t^2 + 98 t + 1), $$
$$ -1 + 255 t + 3615 t^2  + 8885 t^3 - 255 t^4  + t^5 , $$
$$(t - 1) (t^4  + 121 t^3  - 6494 t^2  + 121 t + 1),$$
$$-1 - 495 t + 4365 t^2  - 6615 t^3  - 380 t^4  + t^5 ,$$
$$ (t - 1)^5 $$

\begin{tikzpicture}
  \node at (0,0) (O) {$\bullet$};
\draw [very thin](0,0) circle(5cm);
 
\draw [dotted](0:5)  \foreach \x in {60,120,...,359} {
                -- (\x:5) 
            }-- cycle (-90:5) node[below=1cm] {$n=6$} ;
\draw [->>,thick](180:6) -- (0:6)  (0:5) node[above=5mm] {$x_0 = x_6$};
\draw [->>,thick](180+60:6) -- (60:6) (60:5) node[above=5mm] {$x_1$};

\draw [->>,thick](180+120:6) -- (120:6) (120:5) node[left=5mm] {$x_2$};
\draw [->>,thick](180+180:6) -- (180:6) (180:5) node[below=5mm] {$x_3$};
\draw [->>,thick](240+180:6) -- (240:6) (240:5) node[below=5mm] {$x_4$};
\draw [->>,thick](300+180:6) -- (300:6) (300:5) node[below=5mm] {$x_5$};
\end{tikzpicture}

$$\mathcal{R}_{{\theta}_{\frac{1}{2}}}: = \left(
\begin{array}{c}
x_1\\
x_2 \\
x_3 \\
x_4 \\
x_5 \\
x_6
\end{array}
\right) \rightarrow \left(
\begin{array}{c}
x_1 \\
x_2  \\
x_3 + 15x_1{\epsilon}^\frac{4}{6} \\
x_4  + 15x_5{\epsilon}^\frac{4}{6}\\
x_5 \\
x_6
\end{array}
\right) \\ \mathcal{R}_{{\theta}_{1}}: = \left(
\begin{array}{c}
x_1\\
x_2 \\
x_3 \\
x_4 \\
x_5 \\
x_6
\end{array}
\right) \rightarrow \left(
\begin{array}{c}
x_1 \\
x_2 \\
x_3 + 6x_2{\epsilon}^\frac{5}{6}  \\
x_4 +  20x_5{\epsilon}^\frac{3}{6} \\
x_5 +  6x_6{\epsilon}^\frac{1}{6} \\
x_6
\end{array}
\right)$$
\begin{tikzpicture}
\draw (0:5) \foreach \x in {51.4286,102.8571,...,359} {
                -- (\x:5)
            }-- cycle (-90:5) node[below=1cm] {$n=7$} ;

\draw [->>,thick](180:6) -- (0:6)  (0:5) node[below=2mm] {$x_0 = x_7$};
\draw [->>,thick](180+51.4286:6) -- (51.4286:6) (51.4286:5) node[below=2mm] {$x_1$};
\draw [->>,thick](180+102.8571:6) -- (102.8571:6) (102.8571:5) node[below=2mm] {$x_2$};
\draw [->>,thick](180+154.2857:6) -- (154.2857:6) (154.2857:5) node[below=2mm] {$x_3$};
\draw [->>,thick](205.7143+180:6) -- (205.7143:6) (205.7143:5) node[below=2mm] {$x_4$};
\draw [->>,thick](-102.8571+180:6) -- (-102.8571:6) (-102.8571:5) node[below=2mm] {$x_5$};
\draw [->>,thick](-51.4286+180:6) -- (-51.4286:6) (-51.4286:5) node[below=2mm] {$x_6$};

\end{tikzpicture}

$$\mathcal{R}_{{\theta}_{\frac{1}{2}}}: = \left(
\begin{array}{c}
x_1\\
x_2 \\
x_3 \\
x_4 \\
x_5 \\
x_6 \\
x_7
\end{array}
\right) \rightarrow \left(
\begin{array}{c}
x_1 \\
x_2 \\
x_3 + 21x_1{\epsilon}^\frac{5}{7} \\
x_4 + 35x_7{\epsilon}^\frac{3}{7} \\
x_5 + 7x_6{\epsilon}^\frac{1}{7}\\
x_6 \\
x_7 
\end{array}
\right)\\ \mathcal{R}_{{\theta}_{1}}: = \left(
\begin{array}{c}
x_1\\
x_2 \\
x_3 \\
x_4 \\
x_5 \\
x_6 \\
x_7
\end{array}
\right) \rightarrow \left(
\begin{array}{c}
x_1 \\
x_2 \\
x_3 + 7x_2{\epsilon}^\frac{-6}{7} \\
x_4 + 35x_1{\epsilon}^\frac{4}{7} \\
x_5 + 21x_7{\epsilon}^\frac{2}{7} \\
x_6 \\
x_7 
\end{array}
\right) $$

$$ \mathcal{M}_{\frac{1}{2}}: = \left(
\begin{array}{ccccccc}
1&0&21{\epsilon}^{\frac{5}{7}}&-112{\epsilon}^{\frac{4}{7}}&378{\epsilon}^{\frac{3}{7}}&1638{\epsilon}^{\frac{2}{7}}&-728{\epsilon}^{\frac{1}{7}} \\
7{\epsilon}^{\frac{1}{7}}&1&7{\epsilon}^{\frac{6}{7}}&-28{\epsilon}^{\frac{5}{7}}&84{\epsilon}^{\frac{4}{7}}&378{\epsilon}^{\frac{3}{7}}&-224{\epsilon}^{\frac{2}{7}}\\
0&0&1&7{\epsilon}^{\frac{6}{7}}&-28{\epsilon}^{\frac{5}{7}}&-112{\epsilon}^{\frac{4}{7}}&35{\epsilon}^{\frac{3}{7}}\\
0&0&0&1&7{\epsilon}^{\frac{6}{7}}&21{\epsilon}^{\frac{5}{7}}&0 \\
0&0&0&0&1&0&0 \\
0&0&0&0&7{\epsilon}^{\frac{1}{7}}&1&0 \\
0&0&0&35{\epsilon}^{\frac{3}{7}}&-224{\epsilon}^{\frac{2}{7}}&-728{\epsilon}^{\frac{1}{7}}&1
\end{array}
\right)$$

$$ \mathcal{M}_{1}: = \left(
\begin{array}{ccccccc}
-4752&-9009{\epsilon}^{\frac{6}{7}}&52745{\epsilon}^{\frac{5}{7}}&-18634{\epsilon}^{\frac{4}{7}}&378{\epsilon}^{\frac{3}{7}}&-1008{\epsilon}^{\frac{2}{7}}&-2310{\epsilon}^{\frac{1}{7}} \\
-924{\epsilon}^{\frac{1}{7}}&1716&10395{\epsilon}^{\frac{6}{7}}&-3850{\epsilon}^{\frac{5}{7}}&84{\epsilon}^{\frac{4}{7}}&-210{\epsilon}^{\frac{3}{7}}&462{\epsilon}^{\frac{2}{7}}\\
462{\epsilon}^{\frac{2}{7}}&-924{\epsilon}^{\frac{1}{7}}&4950&1575{\epsilon}^{\frac{6}{7}}&-28{\epsilon}^{\frac{5}{7}}&84{\epsilon}^{\frac{4}{7}}&-210{\epsilon}^{\frac{3}{7}}\\
-210{\epsilon}^{\frac{3}{7}}&462{\epsilon}^{\frac{2}{7}}&-2100{\epsilon}^{\frac{1}{7}}&540&7{\epsilon}^{\frac{6}{7}}&-28{\epsilon}^{\frac{5}{7}}&84{\epsilon}^{\frac{4}{7}}\\
84{\epsilon}^{\frac{4}{7}}&-210{\epsilon}^{\frac{3}{7}}&756{\epsilon}^{\frac{2}{7}}&-140{\epsilon}^{\frac{1}{7}}&1&7{\epsilon}^{\frac{6}{7}}&-28{\epsilon}^{\frac{5}{7}} \\
560{\epsilon}^{\frac{5}{7}}&-1386{\epsilon}^{\frac{4}{7}}&5082{\epsilon}^{\frac{3}{7}}&-959{\epsilon}^{\frac{2}{7}}&7{\epsilon}^{\frac{1}{7}}&-48&-189{\epsilon}^{\frac{6}{7}} \\
-5775{\epsilon}^{\frac{6}{7}}&12320{\epsilon}^{\frac{5}{7}}&-59059{\epsilon}^{\frac{4}{7}}&316107{\epsilon}^{\frac{3}{7}}&-224{\epsilon}^{\frac{2}{7}}&840{\epsilon}^{\frac{1}{7}}&-2400
\end{array}
\right)$$

The characteristic polynomials for the intermediate matrices again for $\epsilon = -1$ are calculated:\\

$$(t - 1)^7,(t - 1)^7,(t - 1)^7,(t - 1)^7,(t - 1)^7, (t - 1)^7 $$ 
$$ (t - 1)^7,(t - 1) (t^2 + 439 t + 1) (t^2  + 1223 t + 1) (t^2  + 47 t + 1) $$ 
$$ t^7  - 6867 t^6 + 5881099 t^5  + 104325473 t^4  + 371279699 t^3  - 3831674 t^2  + 7210 t - 1 $$
$$ (t - 1)(t^6  + 11656 t^5  + 6716984 t^4  - 278167532 t^3  + 6716984 t^2  + 11656 t^7 + 1)$$ 
$$ t^7- 6524 t^6  - 961751 t^5  - 147862190 t^4  - 137831085 t^3  - 6530741 t^2 + 10983 t - 1, $$ 
$$ (t - 1) (t^6  + 1709 t^5  - 5178599 t^4  - 2234665 t^3  - 5178599 t^2  + 1709 t + 1),$$ 
$$t^7  + 14399 t^6  - 7670488 t^5  - 4987255 t^4  - 126708624 t^3  + 5101761 t^2+ 12698 t - 1,$$ 
$$(t - 1)^7$$

A different approach which leads to the following set of matrices:
$$ \mathcal{M}: = \left(
\begin{array}{ccccccc}
-48&-7&-959&-5082&-1386&-560&-189\\
7&1&140&756&210&84&28\\
28&7&540&2100&462&210&84 \\
84&28&1575&4950&924&462&210 \\
210&84&3850&10395&1716&924&462 \\
-1008&-378&-18634&-52745&69009&-4752&-2310\\
-84&-224&-16107&-59059&-12320&-5775&-2400 
\end{array}
\right)$$
The left normalising matrix $\mathcal{L}$ and the right normalising matrix $\mathcal{R}$
$$ \mathcal{L}: = \left(
\begin{array}{ccccccc}
 1&0&0&0&0&0&0 \\
-1&1&0&0&0&0&0 \\
1&-2&1&0&0&0&0 \\
-1&3&-3&1&0&0&0 \\
1&-4&6&3&1&0&0 \\
-1&5&25&3&2&1&0 \\
1&15&15&1&1&1&1
\end{array}
\right)    \mathcal{R}: = \left(
\begin{array}{ccccccc}
  1&0&0&0&0&0&0 \\
-1&1&0&0&0&0&0 \\
1&-2&1&0&0&0&0 \\
-1&3&4&1&0&0&0 \\
20&40&4&3&1&0&0 \\
-1&5&25&3&2&1&0 \\
15&15&1&1&1&1&1
\end{array}
\right) $$
There are relations between the left and right normalising matrices:
multiplication of the inverse of the right normalising matrix with the left normalising matrix sandwiched in between by $J(p,0)+J(p,1)$ and $P(p,-1)$ gives identity, where $J(p,0)$ is the identity matrix of order $p$ , $J(p,r)$ is the $r$-shifted upper diagonal matrix and $P(p,r)$ is the $r$-shifted matrix.

\begin{tikzpicture}
\draw (0:5) \foreach \x in {45,90,...,359} {
                -- (\x:5)
            } -- cycle (-90:5) node[below = 1cm] {$n=8$} ;
\draw [->>,thick](180:6) -- (0:6) (0:5) node[below=2mm] {$x_0 = x_8$}; 
\draw [->>,thick](180+45:6) -- (45:6) (45:5) node[above=2mm] {$x_1$}; 
\draw [->>,thick](90+180:6) -- (90:6) (90:5) node[above=2mm] {$x_2$};
\draw [->>,thick](135+180:6) -- (135:6) (135:5) node[above=2mm] {$x_3$};
\draw [->>,thick](180+180:6) -- (180:6) (180:5) node[below=2mm] {$x_4$};
\draw [->>,thick](225+180:6) -- (225:6) (225:5) node[below=2mm] {$x_5$};
\draw [->>,thick](180+270:6) -- (270:6) (270:5) node[below=2mm] {$x_6$}; 
\draw [->>,thick](315+180:6) -- (315:6) (315:5) node[below=2mm] {$x_7$};

\end{tikzpicture}

$$\mathcal{R}_{{\theta}_{\frac{1}{2}}}: = \left(
\begin{array}{c}
x_1\\
x_2 \\
x_3 \\
x_4 \\
x_5 \\
x_6 \\
x_7 \\
x_8
\end{array}
\right) \rightarrow \left(
\begin{array}{c}
x_1 \\
x_2 \\
x_3 + 8x_2{\epsilon}^\frac{7}{8}  \\
x_4 + 56x_1{\epsilon}^\frac{5}{8} \\
x_5 + 56x_8{\epsilon}^\frac{3}{8} \\
x_6 + 8x_7{\epsilon}^\frac{1}{8} \\
x_7 \\
x_8
\end{array}
\right) \\ \mathcal{R}_{{\theta}_{1}}: = \left(
\begin{array}{c}
x_1\\
x_2 \\
x_3 \\
x_4 \\
x_5 \\
x_6 \\
x_7 \\
x_8
\end{array}
\right) \rightarrow \left(
\begin{array}{c}
x_1 \\
x_2 \\
x_3 \\
x_4 + 28x_2{\epsilon}^\frac{3}{4} \\
x_5 + 70x_1{\epsilon}^\frac{1}{2} \\
x_6 + 28x_8{\epsilon}^\frac{1}{4} \\
x_7 \\
x_8
\end{array}
\right)$$


\begin{tikzpicture}
\draw (0:5) \foreach \x in {40,80,...,360} {
                -- (\x:5)
            } -- cycle (-90:5) node[below = 1cm] {$n=9$} ;
\draw [->>,thick](180:6)--(0:6) (0:5) node[below=2mm] {$x_0 = x_9$};
\draw [->>,thick](180+40:6) -- (40:6) (40:5) node[above=2mm] {$x_1$};
\draw [->>,thick](80+180:6) -- (80:6) (80:5) node[above=2mm] {$x_2$};
\draw [->>,thick](120+180:6) -- (120:6) (120:5) node[above=2mm] {$x_3$};
\draw [->>,thick](160+180:6) -- (160:6) (160:5) node[left=2mm] {$x_4$};
\draw [->>,thick](200+180:6) -- (200:6) (200:5) node[below=2mm] {$x_5$};
\draw [->>,thick](240+180:6) -- (240:6) (240:5) node[below=2mm] {$x_6$};
\draw [->>,thick](280+180:6) -- (280:6) (280:5) node[below=2mm] {$x_7$};
\draw [->>,thick](320+180:6) -- (320:6) (320:5) node[below=2mm] {$x_8$};
\end{tikzpicture}

$$\mathcal{R}_{{\theta}_{\frac{1}{2}}}: = \left(
\begin{array}{c}
x_1\\
x_2 \\
x_3 \\
x_4 \\
x_5 \\
x_6 \\
x_7 \\
x_8 \\
x_9
\end{array}
\right) \rightarrow \left(
\begin{array}{c}
x_1 \\
x_2 \\
x_3 + 9x_2{\epsilon}^\frac{8}{9}  \\
x_4 + 84x_1{\epsilon}^\frac{2}{3} \\
x_5 + 126x_9{\epsilon}^\frac{4}{9} \\
x_6 + 36x_8{\epsilon}^\frac{2}{9}  \\
x_7 \\
x_8 \\
x_9 
\end{array}
\right) \mathcal{R}_{{\theta}_{1}}: = \left(
\begin{array}{c}
x_1\\
x_2 \\
x_3 \\
x_4 \\
x_5 \\
x_6 \\
x_7 \\
x_8 \\
x_9
\end{array}
\right) \rightarrow \left(
\begin{array}{c}
x_1 \\
x_2 \\
x_3 \\
x_4 + 36x_2{\epsilon}^\frac{7}{9}  \\
x_5 + 126x_1{\epsilon}^\frac{5}{9} \\
x_6 + 84x_9{\epsilon}^\frac{1}{3}  \\
x_7 + 9x_8{\epsilon}^\frac{1}{9} \\
x_8 \\
x_9 
\end{array}
\right)$$

The form of the matrix for half$\mathcal{M}_{\frac{1}{2}}$ and full turn $\mathcal{M}_{1}$ for $\epsilon = -1$ are given as follows: (we denote $-1$as $\epsilon$ everywhere:)

$$ \mathcal{M}_{\frac{1}{2}}: = \left(
\begin{array}{ccccccccc}
 1&0&0&84{\epsilon}^{\frac{2}{3}}&-630{\epsilon}^{\frac{5}{9}}&2772{\epsilon}^{\frac{4}{9}}&-9240{\epsilon}^{\frac{1}{3}}&42372{\epsilon}^{\frac{2}{9}}&10575{\epsilon}^{\frac{4}{9}}\\
-315{\epsilon}^{\frac{1}{9}}&1&9{\epsilon}^{\frac{8}{9}}&-45{\epsilon}^{\frac{7}{9}}&165{\epsilon}^{\frac{2}{3}}&-495{\epsilon}^{\frac{5}{9}}&1287{\epsilon}^{\frac{4}{9}}&-9240{\epsilon}^{\frac{1}{3}}&4950{\epsilon}^{\frac{2}{9}}\\
36{\epsilon}^{\frac{2}{9}}&0&1&9{\epsilon}^{\frac{8}{9}}&-45{\epsilon}^{\frac{7}{9}}&165{\epsilon}^{\frac{2}{3}}&-495{\epsilon}^{\frac{5}{9}}&2772{\epsilon}^{\frac{4}{9}}&-1050{\epsilon}^{\frac{1}{3}}\\
0&0&0&1&9{\epsilon}^{\frac{8}{9}}&-45{\epsilon}^{\frac{7}{9}}&165{\epsilon}^{\frac{2}{3}}&-630{\epsilon}^{\frac{5}{9}}&126{\epsilon}^{\frac{4}{9}}\\
0&0&0&0&1&9{\epsilon}^{\frac{8}{9}}&-45{\epsilon}^{\frac{7}{9}}&84{\epsilon}^{\frac{2}{3}}&0 \\
0&0&0&0&0&1&9{\epsilon}^{\frac{8}{9}}&0&0 \\
0&0&0&0&0&0&1&0&0 \\
0&0&0&0&0&36{\epsilon}^{\frac{2}{9}}&-315{\epsilon}^{\frac{1}{9}}&1&0 \\
0&0&0&0&126{\epsilon}^{\frac{4}{9}}&-1050{\epsilon}^{\frac{1}{3}}&4950{\epsilon}^{\frac{2}{9}}&-10575{\epsilon}^{\frac{1}{9}}&1 \\
\end{array}
\right)
$$
$\mathcal{M}_{1}:=$
\\
$$\left(
\begin{array}{ccccccccc}
 140140&288288{\epsilon}^{\frac{8}{9}}&{\epsilon}^{\frac{7}{9}}&-5384652{\epsilon}^{\frac{2}{3}}&1853838{\epsilon}^{\frac{5}{9}}&-80388{\epsilon}^{\frac{4}{9}}&-9240{\epsilon}^{\frac{1}{3}}&25740{\epsilon}^{\frac{2}{9}}&-63063{\epsilon}^{\frac{1}{9}}\\
-12870{\epsilon}^{\frac{1}{9}}&24310&-175032{\epsilon}^{\frac{8}{9}}&524160{\epsilon}^{\frac{7}{9}}&-210210{\epsilon}^{\frac{2}{3}}&11088{\epsilon}^{\frac{5}{9}}&1287{\epsilon}^{\frac{4}{9}}&-3003{\epsilon}^{\frac{1}{3}}&6435{\epsilon}^{\frac{2}{9}}\\
6435{\epsilon}^{\frac{2}{9}}&-12870{\epsilon}^{\frac{1}{9}}&-91520&-252252{\epsilon}^{\frac{8}{9}}&-91728{\epsilon}^{\frac{7}{9}}&-4290{\epsilon}^{\frac{2}{3}}&-495{\epsilon}^{\frac{5}{9}}&1287{\epsilon}^{\frac{4}{9}}&-3003{\epsilon}^{\frac{1}{3}}\\
-3003{\epsilon}^{\frac{1}{3}}&6435{\epsilon}^{\frac{2}{9}}&45045{\epsilon}^{\frac{1}{9}}&-112112&-36036{\epsilon}^{\frac{8}{9}}&1440{\epsilon}^{\frac{7}{9}}&165{\epsilon}^{\frac{2}{3}}&-495{\epsilon}^{\frac{5}{9}}&1287{\epsilon}^{\frac{4}{9}}\\
1287{\epsilon}^{\frac{4}{9}}&-3003{\epsilon}^{\frac{1}{3}}&-20592{\epsilon}^{\frac{2}{9}}&45045{\epsilon}^{\frac{1}{9}}&-12320&-396{\epsilon}^{\frac{8}{9}}&-45{\epsilon}^{\frac{7}{9}}&165{\epsilon}^{\frac{2}{3}}&-495{\epsilon}^{\frac{5}{9}} \\
-495{\epsilon}^{\frac{5}{9}}&1287{\epsilon}^{\frac{4}{9}}&8580{\epsilon}^{\frac{1}{3}}&-15840{\epsilon}^{\frac{2}{9}}&3465{\epsilon}^{\frac{1}{9}}&-80&9{\epsilon}^{\frac{8}{9}}&-45(-1)^{\frac{7}{9}}&165{\epsilon}^{\frac{2}{3}} \\
165{\epsilon}^{\frac{2}{3}}&-495{\epsilon}^{\frac{5}{9}}&-3168{\epsilon}^{\frac{4}{9}}&4620{\epsilon}^{\frac{1}{3}}&-720{\epsilon}^{\frac{2}{9}}&9{\epsilon}^{\frac{1}{9}}&1&9{\epsilon}^{\frac{8}{9}}&-45{\epsilon}^{\frac{7}{9}} \\
-16380{\epsilon}^{\frac{7}{9}}&42042{\epsilon}^{\frac{2}{3}}&281358{\epsilon}^{\frac{5}{9}}&-529668{\epsilon}^{\frac{4}{9}}&118344{\epsilon}^{\frac{1}{3}}&-2799{\epsilon}^{\frac{2}{9}}&-315{\epsilon}^{\frac{1}{9}}&1540&5544{\epsilon}^{\frac{8}{9}} \\
126126{\epsilon}^{\frac{8}{9}}&-286650{\epsilon}^{\frac{7}{9}}&-1979250{\epsilon}^{\frac{2}{3}}&4502484{\epsilon}^{\frac{5}{9}}&-1286424{\epsilon}^{\frac{4}{9}}&43500{\epsilon}^{\frac{1}{3}}&4950{\epsilon}^{\frac{2}{9}}&-17325{\epsilon}^{\frac{1}{9}}&50050 \\
\end{array}
\right)
$$

The characteristic polynomials of the intermediate matrices (for $\epsilon = -1 $)are as follows:

$$(t - 1)^9 , (t - 1)^9 , (t - 1)^9 , $$
$$(t - 1)^9 , (t - 1)^9 , (t - 1)^9 , $$
$$(t - 1)^9 ,(t - 1)^9 ,(t - 1)^9 , $$
$$ (t - 1) (t^2  + 1294 t + 1)(t^2  + 79 t + 1) (t^2  + 15874 t + 1)(t^2  + 7054 t + 1),$$
$$ t^9  - 147735 t^8  + 2970065547 t^7  - 3170023205223 t^6+675828329399154 t^5  $$
$$ + 147081081042558 t^4  + 6741537707415 t^3 - 2463421239 t^2  + 148707 t - 1, $$
$$ (t - 1) (t^8  + 393679 t^7  + 11615124304 t^6+ 19299925922101 t^5  $$
$$- 1152028241052050 t^4  + 19299925922101 t^3+ 11615124304 t^2  + 393679 t + 1),$$
$$ t^9  - 602397 t^8  + 14430588156 t^7- 32570677137612 t^6  - 501642805172094 t^5 $$
$$ - 889035677638326 t^4 - 10413123392220 t^3  - 16402547844 t^2  + 561789 t - 1, $$
$$(t - 1) (t^8 + 500347 t^7  + 6035901139 t^6  - 5067510858641 t^5 $$
$$ + 161677286153188 t^4 - 5067510858641 t^3  + 6035901139 t^2  + 500347 t + 1), $$
$$ t^9  - 383463 t^8 + 3175966269 t^7  + 1760986951161 t^6 + 55212303757689 t^5 $$
$$ - 44378330827005 t^4  + 3651135711771 t^3  - 2896114941 t^2  + 64863 t - 1,$$
$$(t - 1) (t^8  + 24301 t^7  - 3978173600 t^6  - 1121174185736 t^5 + 31682357496619 t^4 $$
$$ - 1121174185736 t^3  - 3978173600 t^2  + 24301 t + 1), $$
$$ t^9 - 407772 t^8  - 5635914723 t^7  + 2577834372561 t^6  + 75567247178913 t^5 $$
$$+ 34683494824953 t^4  + 3563205077937 t^3  + 7515569187 t^2  - 432063 t - 1, $$
$$  (t - 1)^9 $$


%% file: ode_8.tex
\subsection{T and $\gamma$ polynomials}

For complex values of $p$ that is denoted by $\alpha$ in what follows, the \textit{nir}-transform satisfies ODEs of infinite order. We take into consideration the Proposition 6.1 of [ES].

\noindent \textbf{T-polynomials}
\begin{eqnarray*}
T_{0}(\beta)&=& 1
\\
T_{2}(\beta)&=&(1/12)\,\beta
\\
T_{4}(\beta)&=&(1/240)\,\beta\,(-2+5\,\beta)
\\
T_{6}(\beta)&=&(1/4032)\,\beta\,(16+42\,\beta+35\,\beta^{2})
\\
T_{8}(\beta)&=&(1/34560)\,\beta\,(-4+5\,\beta)(36-56\,\beta+35\,\beta^{2})
\\
T_{10}(\beta)&=&(1/101376)\,\beta\, (768-2288\,\beta+2684\,\beta^{2}-1540\,\beta^{3}+385\,\beta^{4})
\\
T_{12}(\beta)&=&(1/50319360)\,\beta\,(-1061376+3327584\,\beta-4252248\,\beta^{2}+2862860\,\beta^{3}
\\
&&-1051050\,\beta^{4}+175175\,\beta^5) 
\\
T_{14}(\beta)&=&(1/6635520)\,\beta\,(552960-1810176\,\beta+2471456\,\beta^{2}-1849848\,\beta^{3}
\\
&& +820820\,\beta^{4}-210210\,\beta^{5}+25025\,\beta^{6}) 
\\
T_{16}(\beta)&=&(1/451215360)\,\beta\,(-200005632+679395072\,\beta-978649472\,\beta^{2}+792548432\,\beta^{3}
\\
&&-397517120\,\beta^{4}+ 125925800\,\beta^{5}-23823800\,\beta^{6}+2127125\,\beta^7)
\\
T_{18}(\beta)&=&(1/42361159680)\,\beta\,(129369047040-453757851648\,\beta+683526873856\,\beta^{2}
\\
&&-589153364352\,\beta^{3}+323159810064\,\beta^{4}-117327450240\,\beta^{5}
\\
&&+27973905960\,\beta^{6}-4073869800\,\beta^{7}+282907625\,\beta^{8})
\\ 
T_{20}(\beta)&=& (1/1471492915200)\,\beta\,(-38930128699392+140441050828800\,\beta-219792161825280\,\beta^{2}
\\
&&+199416835425280\,\beta^{3}-117302530691808\,\beta^{4}+47005085727600\,\beta^{5}-12995644662000\,\beta^{6}
\\
&&+2422012593000\,\beta^{7}-280078548750\,\beta^{8}+15559919375\,\beta^{9})
\\ 
T_{22}(\beta)&=& (1/1758147379200)\,\beta\,(494848416153600-1830317979303936\,\beta+2961137042841600\,\beta^{2}
\\
&&-2805729689044480\,\beta^{3}+1747214980192000\,\beta^{4}-755817391389984\,\beta^{5}
\\
&& +232489541684400\,\beta^{6}-5074916606600\,\beta^{7}+7607466867000\,\beta^{8}
\\
&& -715756291250\,\beta^{9}+32534376875\,\beta^{10})
\\
\end{eqnarray*}

\newpage
\begin{eqnarray*}
T_{24}(\beta)&=&(1/417368899584000)\,\beta\,(-1505662706987827200+5695207005856038912\,\beta
\\
&&-9487372599204065280\,\beta^{2}+9332354263294766080\,\beta^{3}-6096633539052376320\,\beta^{4}
\\
&& +2806128331871953088\,\beta^{5}-937291839756592320\,\beta^{6}+229239926321406000\,\beta^{7}
\\
&& -40598842049766000\,\beta^{8}+5005999501002500\,\beta^{9}-390802935022500\,\beta^{10}
\\
&& +14803141478125\,\beta^{11})
\\ 
T_{26}(\beta)&=& (1/15410543984640)\,\beta\,(844922884529848320-3261358271400247296\,\beta
\\
&& +5576528334428209152\,\beta^{2}-5668465199488266240\,\beta^{3}+3858582205451484160\,\beta^{4}
\\
&& -1870620248833400064\,\beta^{5}+667822651436228288\,\beta^{6}-178292330746770240\,\beta^{7}
\\
&&+35600276746834800\,\beta^{8}-5225593531158000\,\beta^{9}+539680243602500\,\beta^{10}
\\
&& -35527539547500\,\beta^{11}+1138703190625\,\beta^{12}) 
\\
T_{28}(\beta)&=&(1/141874849382400)\,\beta\,(-138319015041155727360+543855095595477762048\,\beta
\\
&& -952027796641042464768\,\beta^{2}+996352286992030556160\,\beta^{3}-703040965960031795200\,\beta^{4}
\\
&& +356312537387839432192\,\beta^{5}-134466795172062184832\,\beta^{6}+38526945410311117760\,\beta^{7}
\\
&& -8436987713444690400\,\beta^{8}+1404048942958662000\,\beta^{9}-173777038440005000\,\beta^{10}
\\
&& +15258232341852500\,\beta^{11}-858582205731250\,\beta^{12}+23587423234375\,\beta^{13})
\\ 
T_{30}(\beta)&=& (1/28026642660065280)\,\beta\,(562009739464769840087040-
\\
&& 2247511941596311764074496\,\beta+4019108379306905439830016\,\beta^{2}
\\
&& -4317745925208072594259968\,\beta^{3}+3145163776677939429416960\,\beta^{4}
\\
&& -1656917203539032341530624\,\beta^{5}+655643919364420586023424\,\beta^{6}
\\
&& -199227919419039256217472\,\beta^{7}+46995751664475880185920\,\beta^{8}
\\
&& -8614026107092938211680\,\beta^{9}+1214778349162323946000\,\beta^{10}
\\
&& -128587452922193265000\,\beta^{11}+9720180867524627500\,\beta^{12}
\\
&& -472946705787806250\,\beta^{13}+11260635852090625\,\beta^{14})
\\ 
T_{32}(\beta)&=& (1/922166952040857600)\,\beta\,(-435617713657079143134658560+
\\
&& 1769501101795425588791476224\,\beta-3226607738514112740810817536\,\beta^2
\\
&& +3549464948895368918909386752\,\beta^3 -2660240772214473252783390720\,\beta^4
\\
&& +1450197119150695223028064256\,\beta^5-597982476958827532820312064\,\beta^6
\\
&& +191042987397310700759115008\,\beta^7-47937412602413871076700160\,\beta^8
\\
&& +9496768069402188916140800\,\beta^9-1480818668624216662540800\,\beta^{10}
\\
&& +179407801092762942700000\,\beta^{11}-16436708736371283360000\,\beta^{12}
\\
&& +1081201211974333450000\,\beta^{13}-45943394276529750000\,\beta^{14}
\\
&& +957154047427703125\,\beta^{15}) 
\\
\end{eqnarray*}
\newpage

\begin{eqnarray*}
T_{34}(\beta)&=& (1/19725496300339200)\,\beta\,(249245942713958501292441600
\\
&& -1027227011317383671945625600\,\beta+1907002180538686183638564864\,\beta^2
\\
&& -2143650324846265987318677504\,\beta^3+1648507262010394952489828352\,\beta^4
\\
&& -926543387509564170390077440\,\beta^5+396194506712421055863955456\,\beta^6
\\
&& -132202422240130962057834496\,\beta^7+34964305711242759446603008\,\beta^8
\\
&& -7388542677376642277816320\,\beta^9+1249066025676468611772160\,\beta^{10}
\\
&& -167917708706445592294400\,\beta^{11}+17680052341202751500000\,\beta^{12}
\\
&& -1416281700897823760000\,\beta^{13}+81822997425819650000\,\beta^{14}
\\
&& -3062892951768650000\,\beta^{15}+56303179260453125\,\beta^{16}) 
\\
T_{36}(\beta) &=& (1/2163255728265599385600)\,\beta\,(-823939351844726605927049141944320
\\
&& +3441870155335118292992998571507712\,\beta-6496555740332629910708455118733312\,\beta^2
\\
&& +7449087490112312991527181748273152\,\beta^3-5864346479278651486249632960479232\,\beta^4
\\
&& +3388164372699248551771666420400128\,\beta^5-1496540121213913485448632812126208\,\beta^6
\\
&& +518865633562168075157007576711168\,\beta^7-143629918822573878550516045152768\,\beta^8
\\
&& +32064170141779614103830545460480\,\beta^9-5796758938020445202118146035200\,\beta^{10}
\\
&& +847300686907144962577480595200\,\beta^{11}-99312783205522355459126136000\,\beta^{12}
\\
&& +9177047640834925292701260000\,\beta^{13}-648457545044826396980100000\,\beta^{14}
\\
&& +33170257743163225434750000\,\beta^{15}-1102176285769664136281250\,\beta^{16}
\\
&& +18009416434144838828125\,\beta^{17}) 
\\
T_{38}(\beta)&=&(1/36926129074234982400)\,\beta\,(474532164251036578163360740147200-
\\
&& 2007451486703515474962315120476160\,\beta+3847934197291457612390714031734784\,\beta^2
\\
&& -4493689786808060863818452674019328\,\beta^3
\\
&& +3614524930468188702037700647059456\,\beta^4-2141316875360841710557998512996352\,\beta^5
\\
&& +973848254246503892750214050676736\,\beta^6-349368584100039649120827420721152\,\beta^7
\\
&& +100669622071229591817492908396544\,\beta^8-23568253290977786453992670028288\,\beta^9
\\
&& +4510878698049367240357658177280\,\beta^{10}-706768823732042351913226513920\,\beta^{11}
\\
&& +90309415620240337258651513600\,\beta^{12}-9317686151162890530973944000\,\beta^{13}
\\
&& +762077701556459749521060000\,\beta^{14}-47865941553425413918500000\,\beta^{15}
\\
&& +2183358737524668003300000\,\beta^{16}-64833899162921419781250\,\beta^{17}
\\
&& +947864022849728359375\,\beta^{18}) 
\end{eqnarray*}
\newpage

\noindent \textbf{The $\gamma$ polynomials}
\begin{eqnarray*}
\gamma_{0}({\alpha}^2) &=& 1 
\\
{\gamma}_{1}({\alpha}^2) &=& (\alpha -1)\,(\alpha+1)
\\
{\gamma}_{2}({\alpha}^2) &=& (\alpha -1)\,(\alpha +1)\,({\alpha}^{2}+23) 
\\
\gamma_{3}({\alpha}^2) &=& (1/5)\,(\alpha-1)\,(\alpha+1)\,(5\,{\alpha}^{4}-298\,\alpha^{2}+11237) 
\\
\gamma_{4}({\alpha}^2)&=& (1/5)\,(\alpha-1)\,(\alpha+1)\,(5\,{\alpha}^{6}-1887\,{\alpha}^{4}-241041\,{\alpha}^{2}+2482411) 
\\ 
{\gamma}_{5}({\alpha}^2) &=& (1/7)\,(\alpha-1)\,(\alpha+1)\,(7\,{\alpha}^{8}-7420\,{\alpha}^{6}
\\
&& +1451274\,{\alpha}^{4}-220083004\,{\alpha}^{2}+1363929895)
\\
{\gamma}_{6}({\alpha}^2) &=& (1/35)\,({\alpha}-1)\,({\alpha}+1)\,(35\,{\alpha}^{10}
\\
&& -78295\,{\alpha}^{8}+76299326\,{\alpha}^{6}+25171388146\,{\alpha}^{4}
\\
&& -915974552561\,{\alpha}^{2}+4175309343349) 
\\
{\gamma}_{7}({\alpha}^2) &=& (1/5)\,({\alpha}-1)\,({\alpha}+1)\,(5\,{\alpha}^{12}
\\
&& -20190\,{\alpha}^{10}+45700491\,{\alpha}^{8}-19956117988\,{\alpha}^{6}
\\
&& +7134232164555\,{\alpha}^{4}-142838662997982\,{\alpha}^{2}+525035501918789) 
\\
{\gamma}_{8}({\alpha}^2) &=& (1/5)\,({\alpha}-1)\,({\alpha}+1)\,(5\,{\alpha}^{14}
\\
&& -32963\,{\alpha}^{12}+135959721\,{\alpha}^{10}-286961736847\,{\alpha}^{8}
\\
&& -181220025335249\,{\alpha}^{6}+14569888571515191\,{\alpha}^{4}
\\
&& -201164685264533917\,{\alpha}^{2}+628141988536245979) 
\\
{\gamma}_{9}({\alpha}^2) &=&(1/275)\,({\alpha}-1)\,({\alpha}+1)\,(275\,{\alpha}^{16}
\\
&& -2758360\,{\alpha}^{14}+18264509780\,{\alpha}^{12}-84980828208232,{\alpha}^{10}
\\
&& +66925280510995058\,{\alpha}^{8}-44067073909426812136\,{\alpha}^{6}
\\
&& +1854998841811009566164\,{\alpha}^{4}-19662907822146661450072\,{\alpha}^2
\\
&& +53917386529177385523923) 
\\
{\gamma}_{10}({\alpha}^2) &=& (1/385)\,({\alpha}-1)\,({\alpha}+1)\,(385\,{\alpha}^{18}
\\
&& -5575185\,{\alpha}^{16}+54739917540\,{\alpha}^{14}-441040104588468\,{\alpha}^{12}
\\
&& +1704035085901524414\,{\alpha}^{10}+1756960395376174557138\,{\alpha}^{8}
\\
&& -250051464026304718963308\,{\alpha}^{6}+6980518270250459747738748\,{\alpha}^{4}
\\
&& -60455711535001643407631031\,{\alpha}^{2}+148934765720971351352763767) 
\\
{\gamma}_{11}({\alpha}^2) &=& (1/455)\,({\alpha}-1)\,({\alpha}+1)\,(455\,{\alpha}^{20}
\\
&& -9133670\,{\alpha}^{18}+125862813531\,{\alpha}^{16}
\\
&& -1551851251661256\,{\alpha}^{14}+13138809751628741310\,{\alpha}^{12}
\\
&& -16369514872967259031716\,{\alpha}^{10}+17232759645860853077282430\,{\alpha}^{8}
\\
&&-1244436833719440506308299720\,{\alpha}^{6}+25839020308128868113601951611\,{\alpha}^{4}
\\
&& -190358675897996551624967496038\,{\alpha}^{2}+428338546734334777277256756263)
\\
\end{eqnarray*}
\newpage
\begin{eqnarray*}
{\gamma}_{12}({\alpha}^2) &=&(1/2275)\,({\alpha}-1)\,({\alpha}+1)\,(2275\,{\alpha}^{22}
\\
&& -61286225\,{\alpha}^{20}+1138143234085\,{\alpha}^{18}
\\
&& -19899973184243295\,{\alpha}^{16}+289010332500834015006\,{\alpha}^{14}
\\
&& -1800790721758319530222794\,{\alpha}^{12}-2750182482964852045472807958\,{\alpha}^{10}
\\
&& +609308050185234247965268669122\,{\alpha}^{8}
\\
&& -28432172394701052919004186087217\,{\alpha}^{6}
\\
&& +469816115524203185110937294104043\,{\alpha}^{4}
\\
&& -3027851448645708864699151581524191\,{\alpha}^{2}+6301150244751080741665843707891149) 
\\
{\gamma}_{13}({\alpha}^2) &=& (1/175)\,({\alpha}-1)\,({\alpha}+1)\,(175\,{\alpha}^{24}
\\
&& -6160980\,{\alpha}^{22}+149439244350\,{\alpha}^{20}-3517885582792900\,{\alpha}^{18}
\\
&& +77320157515532388801\,{\alpha}^{16}-1061289073759987813960872\,{\alpha}^{14}
\\
&& +1919804727678589315033522404\,{\alpha}^{12}-2952731541639113644408376751144\,{\alpha}^{10}
\\
&& +324872646459937078198178265809697\,{\alpha}^{8}
\\
&& -11031411694939020431719354111233796\,{\alpha}^{6}
\\
&& +151642180639270311056923944406592190\,{\alpha}^{4}
\\
&& -872621327007244741940501331510316308\,{\alpha}^2
\\
&& +1695881990125518108674571524660426383) 
\\
{\gamma}_{14}({\alpha}^2) &=& (1/25)\,({\alpha}-1)\,({\alpha}+1)\,(25\,{\alpha}^{26}
\\
&& -1125085\,{\alpha}^{24}+34764296190\,{\alpha}^{22}-1061022101650654\,{\alpha}^{20}
\\
&& +32669714654074122547\,{\alpha}^{18}-771636658087676104012503\,{\alpha}^{16}
\\
&& +7084687100305977812202981972\,{\alpha}^{14} 
\\
&& +15014942311637580765119538315852\,{\alpha}^{12}
\\
&& -4768962770078816392189126626103953\,{\alpha}^{10}
\\
&& +332894266374078575368155254033280197\,{\alpha}^{8}
\\
&& -8828958071754076606079100061553173298\,{\alpha}^{6}
\\
&& +104157175058098445988624066193639753842\,{\alpha}^{4}
\\
&& -543550947078313469794032584810082731483\,{\alpha}^{2}
\\
&& +994159279221093204357661985843042030351) 
\end{eqnarray*}

%% file: ODEbis.tex
%
\newpage
\section{ODEs: numerical calculations}

As we have already seen in [ES] and in the course of the last sections, for all polynomial inputs $f$ and all monomial inputs $F$, the inner generators corresponding to them verify ordinary differential equation of linear homogeneous type with polynomial coefficients. In this section we report the negative computational results strongly suggestive of the non-existence of ODEs for functions $F$ that are rational but not monomial, i.e. not of the form ${(1-ax)}^p$, for $p \in \mathbb{Z}$

\subsection{Likely non-existence of ODEs for general rational inputs}

The essential points here are the following: 

The procedure involved in the numerical calculations is based on the exploitation of the \textit{nir}- transform for a given function $f$. We explicitely calculate the \textit{nir}- transform and then subject it to an arbitrary system of differential operators. (here we put the system as dependant on the variable , the function $f$ and the first three derivatives $f^{\prime}$, $f^{\prime\prime}$ and $f^{\prime\prime\prime}$. For as far as the theory already predicts and verified by numerical results too, that the equations are essentially of not such a high degree but of a high order. Care has to be taken about the number of co-efficients to be calculated for the differential system. For the interested reader the program (which is not so complicated) is put at disposal.
 
\begin{itemize}
\item To test the trustworthiness of the method, we first verify, in the case of polynomial inputs $f$, the existence of covariant (resp. variable) ODEs for the semi-entire part (resp. the totality) of the \textit{nir}-transform, and their coincidence with the ODEs predicted by the theory.

\item Then we move on to the existence of ODEs for monomial inputs (i.e. with only on zero or pole but of arbitrary order $p$ 
$$ f(x) := -p\log(1-px)  \,\,\,\,F(x) = (1-px)^p $$ 

or 
$$ f(x) := +p\log(1-px)  \,\,\,\,F(x) = (1+px)^{-p} $$ 

\quad
\begin{eqnarray*}
k(n) &:=& singular\big{(}  \int^{\infty}_{0} e^{-\beta({\partial)}_{\tau}) f(\frac{\tau}{n})}_{\#} d\tau \big{)}    \,\,\,\, \in \Gamma(1/2)n^{1/2}\mathbb{Q}[[n^{-1}]] 
\\
h(\nu) &:=& formal\big{(}  \frac{1}{2 \pi i}\int^{c+i\infty}_{c-i\infty} k(n) e^{\nu n} \frac{dn}{n} \big{)} = h(\nu)   \,\,\,\, \in \nu^{-1/2}\mathbb{Q}\{\nu\} 
\\
\hat{k}(\nu) &:=& formal\big{(} \frac{1}{2 \pi i}\int^{c+i\infty}_{c-i\infty} k(n) e^{\nu n} dn \big{)} = h^{\prime}    \,\,\,\, \in \nu^{-3/2}\mathbb{Q}\{\nu\} 
\end{eqnarray*}

For this case, the global \textit{nir} transforms don't verify any (variable) ODEs. Here only the singular parts which in this case comprises of only half-integral powers verifies covariant ODEs with polynomial coefficients. More precisely , say for $p=1$ we just have a single root of unity (monomial) There is only trivial or elementary resurgence in this situation. The integer part of the  equation is relatively complicated and the numerical tests seem to rule out the existence of at least reasonable simple ODEs. The results have been checked to $N=500$ (number of terms) and for order of the equation $d=5$.

\item $F$ rational but not monomial of the type $\frac{1-x}{1+x}$. This case is the simplest and in particular the only case of a non-monomial, rational $F$ that gives rise to a finite net of singularities. Here, neither the singular nor the regular part of the \textit{nir}-transform seems to verify any ordinary differential equation  
      
\end{itemize}

Explicit numerical computations to this effect have been carried out. In a large domain they negate the existence of any differential equations for the above inputs. To give an idea of the validity of the computations, we give the following technical details:
$d$: the degree (in $\nu$ or $n^{-1}$) of the differential equations that we are looking for. $\delta$ is the order of this differential equation. For our case we have considered the differential equation that would involve at the most the fifth derivative i.e. $\delta = 5$. The number of coefficients $N = (1+\delta)(1+d)$ calculated is 500. 
\\
\begin{center}

 \begin{tabular}{l|c c c r}
 &(2,d) & (3,d) & (4,d) & (5,d)\\
\hline
 N & 500  & 500 & 350 & 300 \\
 n & 30 & 30 & 25 & 20 \\
\end{tabular}
\end{center}

\noindent where $n$ is the number of digits after decimal.

\noindent Remark: The \textit{nir} transform which consists of the nine link chain with the essential nontrivial step that is the so called \textit{mir} transform (and integro-differential operator). A detailed description of this transform is given in [ES].

Remark: For the alternative manner of trying to work out the singularities, we resort directly to the asymptotic Taylor series coefficients. We could always shift the problem to the neighbourhood of the origin, details to appear in [SS].
%


%% file: ode_10.tex
\bigskip

\section{References.}
{\bf [C1]} O. Costin, {\it Asymptotics and Borel summability,}  2008, CRC Press.
{\it  \dots  }
\\
{\bf [C2]} O. Costin, {\it Global reconstruction of analytic functions from local expansions,}
Preprint 2007, Ohio State University.
{\it  \dots  }
\\
{\bf [CG1]} O. Costin and S. Garoufalidis,
{\it Resurgence of the Kontsevich-Zagier power series, }
preprint 2006 math. GT/0609619.
\\
{\bf [CG2]} O. Costin and S. Garoufalidis,
{\it Resurgence of the fractional polylogarithms, }
preprint 2007, math. CA/0701743.
\\
{\bf [CG3]} O. Costin and S. Garoufalidis,
{\it Resurgence of the Euler-MacLaurin summation formula, }
preprint 2007, math. CA/0703641.
\\
{\bf [E1]} J.Ecalle, {\it Les fonctions resurgentes: Vol.1}, Publ. Math. Orsay, 1981\\
{\bf [E2]} J.Ecalle, {\it Les fonctions resurgentes: Vol.2}, Publ. Math. Orsay, 1981\\
{\bf [E3]} J.Ecalle, {\it Les fonctions resurgentes: Vol.3}, Publ. Math. Orsay, 1985\\
{\bf [E4]} J.Ecalle., {\it Six lectures on transseries etc...,}
in: Bifurcations and Periodic Orbits etc , D. Schlomiuk ed., 1993, Kluwer, p 75-184.
\\
{\bf [E5]} J.Ecalle, {\it Recent advances in the analysis of divergence and singularities},
in: Normal Forms etc, Y.Ilyashenko and C.Rousseau eds., 2004, Kluwer,p 87-186.
\\
{\bf [ES]} J.Ecalle and Sh.Sharma, {\it Power series with sum-product Taylor coefficients and
              their resurgence algebra.}
 \\
{\bf [SS]} Sh.Sharma, {\it Some examples of SP-series\,: a theoretical and numerical investigation},
(forthcoming).
\\